\documentclass[a4paper]{article}
\usepackage{amssymb}
\usepackage{amsfonts}
\usepackage{amsmath}
\usepackage{amsthm}
\usepackage{authblk}
\usepackage[USenglish]{babel}
\usepackage{esint}
\usepackage{graphicx}
\usepackage[colorlinks=true]{hyperref}
\usepackage[latin1]{inputenc}
\newtheorem{theorem}{Theorem}[section]
\newtheorem{lemma}[theorem]{Lemma}
\newtheorem{proposition}[theorem]{Proposition}
\newtheorem{corollary}[theorem]{Corollary}
\newtheorem{remark}[theorem]{Remark}
\newtheorem{definition}[theorem]{Definition}

\numberwithin{equation}{section}
\newcommand{\N}{\mathbb N}
\newcommand{\Z}{\mathbb Z}

\newcommand{\R}{\mathbb R}

\renewcommand{\S}{\mathbb S}

\newcommand{\B}{\mathbb B}

\newcommand{\mcal}{\mathcal}

\newcommand{\mrm}{\mathrm}

\renewcommand{\a}{\alpha}
\renewcommand{\b}{\beta}

\renewcommand{\d}{\delta}
\newcommand{\D}{\Delta}
\newcommand{\e}{\varepsilon}
\newcommand{\z}{\zeta}

\newcommand{\la}{\lambda}
\newcommand{\La}{\Lambda}
\newcommand{\s}{\sigma}
\newcommand{\si}{\varsigma}
\newcommand{\Si}{\Sigma}

\newcommand{\ph}{\varphi}

\renewcommand{\O}{\Omega}
\newcommand{\wt}{\widetilde}

\newcommand{\ol}{\overline}

\newcommand{\ub}{\underbrace}
\newcommand{\fr}{\frac}
\newcommand{\pa}{\partial}
\newcommand{\n}{\nabla}
\newcommand{\fa}{\forall}
\newcommand{\ex}{\exists}
\newcommand{\es}{\emptyset}
\newcommand{\wk}{\rightharpoonup}
\newcommand{\inc}{\hookrightarrow}

\newcommand{\us}{\underset}

\newcommand{\To}{\Rightarrow}
\newcommand{\ot}{\leftarrow}

\newcommand{\ltot}{\longleftrightarrow}

\newcommand{\sm}{\setminus}
\renewcommand{\Cup}{\bigcup}

\newcommand{\sub}{\subset}

\newcommand{\sne}{\subsetneq}

\newcommand{\nin}{\not\in}
\newcommand{\eq}{\equiv}
\newcommand{\ox}{\otimes}

\newcommand{\Pl}{\bigoplus}

\newcommand{\x}{\times}
\renewcommand{\c}{\circ}
\newcommand{\cd}{\cdot}
\newcommand{\ds}{\dots}

\newcommand{\tx}{\text}
\newcommand{\q}{\quad}
\renewcommand{\l}{\left}
\renewcommand{\r}{\right}

\newcommand{\bthm}{\begin{theorem}}
\newcommand{\ethm}{\end{theorem}}
\newcommand{\blem}{\begin{lemma}}
\newcommand{\elem}{\end{lemma}}
\newcommand{\bprop}{\begin{proposition}}
\newcommand{\eprop}{\end{proposition}}
\newcommand{\bcor}{\begin{corollary}}
\newcommand{\ecor}{\end{corollary}}
\newcommand{\brem}{\begin{remark}}
\newcommand{\erem}{\end{remark}}
\newcommand{\bdefi}{\begin{definition}}
\newcommand{\edefi}{\end{definition}}
\newcommand{\bpf}{\begin{proof}}
\newcommand{\epf}{\end{proof}}
\newcommand{\bl}{\begin{array}{l}}
\newcommand{\bll}{\begin{array}{ll}}
\newcommand{\barr}{\begin{array}}
\newcommand{\earr}{\end{array}}
\newcommand{\bite}{\begin{itemize}}
\newcommand{\eite}{\end{itemize}}
\newcommand{\bequ}{\begin{equation}}
\newcommand{\eequ}{\end{equation}}
\newcommand{\beqa}{\begin{eqnarray}}
\newcommand{\eeqa}{\end{eqnarray}}
\newcommand{\beqy}{\begin{eqnarray*}}
\newcommand{\eeqy}{\end{eqnarray*}}
\newcommand{\bin}[2]{\left(\genfrac{}{}{0pt}{}{#1}{#2}\right)}
\newcommand{\qm}[1]{``#1''}
\allowdisplaybreaks

\begin{document}

\everymath{\displaystyle}

\title{A general existence result for stationary solutions to the Keller-Segel system}
\author{Luca Battaglia\thanks{Universit\`a degli Studi Roma Tre, Dipartimento di Matematica e Fisica, Largo S. Leonardo Murialdo 1, 00146 Roma - lbattaglia@mat.uniroma3.it}}
\date{}

\maketitle\

\begin{abstract}
\noindent We consider the following Liouville-type PDE, which is related to stationary solutions of the Keller-Segel's model for chemotaxis:
$$\l\{\bll-\D u+\b u=\rho\l(\fr{e^u}{\int_\O e^u}-\fr{1}{|\O|}\r)&\tx{in }\O\\\pa_\nu u=0&\tx{on }\pa\O\\\earr\r.,$$
where $\O\sub\R^2$ is a smooth bounded domain and $\b,\rho$ are real parameters. We prove existence of solutions under some algebraic conditions involving $\b,\rho$. In particular, if $\O$ is not simply connected, then we can find solution for a generic choice of the parameters. We use variational and Morse-theoretical methods.
\end{abstract}\

\section{Introduction}\

We are interesting in the study of the following partial differential equation:

\bequ\label{ks}
\l\{\bll-\D u+\b u=\rho\l(\fr{e^u}{\int_\O e^u}-\fr{1}{|\O|}\r)&\tx{in }\O\\\pa_\nu u=0&\tx{on }\pa\O\\\earr\r.\tag{$P_{\b,\rho}$}.
\eequ

Here, $\O\sub\R^2$ is a smooth bounded open domain in the plane, $\b$ and $\rho$ are real parameters and $|\O|$ is the Lebesgue measure of $\O$.\\

Problem \eqref{ks} is related to a model introduced by Keller and Segel in \cite{ks} to study \emph{chemotaxis} in biology, namely the movement of organisms according to the the presence of chemicals in the environment.\\
In particular, \eqref{ks} models stationary solutions in Keller and Segel's model.\\

In the case $\b>0$, solutions to \eqref{ks} have been found via a mountain-pass argument in \cite{ww}, whereas families of blowing-up solutions have been constructed in \cite{pv,dpv,ap,bcn1,bcn2}.\\
Here we allow the parameter $\b$ to have any sign. We will tackle problem \eqref{ks} variationally; in fact, solutions are all and only the critical points of the energy functional
\bequ\label{j}
\mcal J_{\b,\rho}(u)=\fr{1}2\int_\O\l(|\n u|^2+\b u^2\r)-\rho\log\int_\O e^u.
\eequ
Since both \eqref{ks} and \eqref{j} are invariant under addition of constants, it will not be restrictive to look for solutions to \eqref{ks} satisfying $\int_\O u=0$; equivalently, we will consider $\mcal J_{\b,\rho}$ not on its natural domain $H^1(\O)$, but rather on its subspace
$$\ol H^1(\O):=\l\{u\in H^1(\O):\,\int_\O u=0\r\}.$$\

In particular, we will study the topology and the homology of very low energy sub-levels
$$\mcal J_{\b,\rho}^c:=\l\{u\in\ol H^1(\O):\,\mcal J_{\b,\rho}(u)\le c\r\}$$
with $c=-L\ll0$, then will deduce existence of solutions using Morse theory. This argument has been introduced in \cite{dm} for a fourth-order problem in Riemannian geometry and it has become a rather classical tool in the study of Liouville-type equation (see for instance the surveys \cite{mal2,mal3}).\\
With respect to most previous results, the main new difficulties are given by the Neumann boundary conditions, rather than Dirichlet, and the presence of the extra linear term $\b u$.\\

Neumann conditions may cause concentration on $\pa\O$, which was excluded in the case of Dirichlet conditions, whereas concentration on the interior is similar in the two cases. The main difference when concentration occurs at the boundary is due to the fact that, roughly speaking, a shrinking ball centered at a point $p\in\pa\O$ is asymptotically half of a shrinking ball contained in $\O$.\\
The argument to fix such an issue is inspired by \cite{ndi}, which deals with a similar problem on higher-dimensional manifolds with boundary.\\

Another issue may be given by linear part $-\D+\b$ not being positive definite, if $\b<0$, which is a new feature in second order Liouville equations.\\
This naturally leads to consider the projection of $u$ into positive and negative sub-spaces of the operator $-\D+\b$. Precisely, we take an orthonormal frame $\{\ph_i\}_{i\in\N}$ of eigenfunctions of $-\D$ with associated positive non-decreasing eigenvalues $\{\la_i\}_{i\in\N}$ (counted with multiplicity), so that
$$u=\sum_{i=1}^{+\infty}\l(\int_\O u\ph_i\r)\ph_i\q\q\q\To\q\q\q-\D u=\sum_{i=1}^{+\infty}\la_i\l(\int_\O u\ph_i\r)\ph_i.$$
Now, if $-\b$ is not an eigenvalue of $-\D$, then $-\la_{I+1}<\b<-\la_I$ for some $I\in\N$, therefore we can define the projection $\Pi_I$ on a finite dimensional subspace, on which orthogonal $-\D+\b$ is positive definite:
\bequ\label{pi}
\Pi_Iu:=\l(\int_\O u\ph_1,\ds,\int_\O u\ph_I\r)\in\R^I.
\eequ

Arguing as in \cite{dm}, we can show that, if $\mcal J_{\b,\rho}(u)\ll0$, then either $\fr{e^u}{\int_\O e^u}$ is concentrated around a finite number of points or $\Pi_Iu$ is large (or both occur).\\
To express this alternative we will use the \emph{join}, which has been used in the variational study of Liouville system of two equations, where one has an alternative between the concentration of each component (see \cite{bjmr,bat1,jkm,batmal,bat2}).\\
Given two topological spaces $X$ and $Y$, its join $X\x Y$ is defined as the product between the two spaces and the unit interval, with identifications at each endpoint. We set
\bequ\label{join}
X\star Y:=\fr{X\x Y\x[0,1]}\sim,
\eequ
with $\sim$ being defined by
$$(x,y,0)\sim(x,y',0)\q\fa\,x\in X,\,\fa\,y,y'\in Y,\q\q\q\q\q\q(x,y,1)\sim(x',y,1)\q\fa\,x,x'\in X,\,\fa\,y\in Y.$$
We will suitably choose $X$ and $Y$ as objects to model each alternative. When $t=0$, the whole $Y$ is collapsed, which means only the first alternative occurs, similarly at $t=1$ only $Y$ is left hence we have the second alternative; if $0<t<1$ we see both spaces because both alternatives occur.\\

To be in position to apply such methods, we need some compactness assumptions on the energy functional $\mcal J_{\b,\rho}$.\\
Unfortunately, Palais-Smale condition is not known to hold for problem \eqref{ks}, nor for similar Liouville-type PDEs. Anyway, such problem can be by-passed thanks to a a deformation lemma by \cite{luc} and some compactness of solutions to \eqref{ks} holding locally uniformly in $\rho$.\\
As in most known results, such a result holds true provided $\rho$ is not an integer multiple of $4\pi$. This time we also need $-\b$ not to be an eigenvalue of $-\D$, in order for the linear operator to be non-degenerate (see Section $2$).\\

Finally, we need to verify that the solution found using these tools is not the trivial one $u\eq0$. This is equivalent to evaluate the Morse index of the trivial solution.\\
We will get non-trivial solutions under some algebraic condition involving the parameters $\b,\rho$ and the eigenvalues $\la_i$. In the case when $\O$ is multiply connected, we are allowed to take more cases, since the topology of low sublevels $\mcal J_{\b,\rho}^{-L}$ is more involved.\\
Precisely, the main result of this paper is the following:

\bthm\label{ex}${}$\\
Assume $\b\ne-\la_i\ne\b-\fr{\rho}{|\O|}$ for any $i\in\N$ and $\rho\nin4\pi\N$ and let $I\le J,K$ be non-negative integers such that
$$4K\pi<\rho<4(K+1)\pi\q\q\q-\la_{I+1}<\b<-\la_I\q\q\q-\la_{J+1}<\b-\fr{\rho}{|\O|}<-\la_J.$$
If $\O$ is simply connected and $2K+I\ne J$, then the problem \eqref{ks} has non-trivial solutions.\\
If $\O$ is not simply connected and $(K,I)\ne(0,J)$, then the problem \eqref{ks} has non-trivial solutions.
\ethm\

\brem${}$\\
Since Theorem \ref{ex} is proved via Morse theory, then the same arguments would also give multiplicity of solutions, provided the energy functional $\mcal J_{\b,\rho}$ is a Morse function, as was done with similar problems in \cite{dem,bdm,bat1,jev,bat2,dlr}. However, it is not clear under which conditions on $\O$ this occurs, although we suspect it is a somehow generic conditions.\\
One easily see that, assuming $u\eq0$ to be the only solution to \eqref{ks}, $\mcal J_{\b,\rho}$ is Morse if and only if $\beta$ and $\rho$ satisfy an algebraic relation (see Proposition \ref{morse}); this fact implies that the solution found in Theorem \ref{ex} is not trivial.\\

The same arguments as Theorem \ref{ex} are also useful, with minor modifications, to find a solution to
$$\l\{\bll-\D u+\b u=\rho\l(\fr{he^u}{\int_\Si he^u}-\fr{1}{|\Si|}\r)&\tx{in }\Si\\\pa_\nu u=0&\tx{on }\pa\Si\\\earr\r.,$$
with $\Si$ being a compact surface with boundary and $h\in C^\infty(\Si)$ being strictly positive.\\
Here, if $h$ is not constant, we can also cover the case in which Theorem \ref{ex} gave a trivial solution. Moreover, by arguing as in the previous references, one can show that $\mcal J_{\b,\rho}$ satisfies the Morse property under a generic assumption on $h$ and/or on the metric $g$ on $\Si$.
\erem\

The plan of this paper is the following. Section $2$ is devoted to the study of compactness properties of the equation \eqref{ks}; Section $3$, which is divided in three sub-sections, deals with the analysis of energy sublevels of $\mcal J_{\b,\rho}$, and finally in Section $4$ Theorem \eqref{ex} is proved.\\

We will denote as $d(x,y)$ the distance between two points $x,y\in\ol\O$ and, similarly, for $\O',\O''\sub\O$,
$$d(x,\O'):=\inf\{d(x,y):\,y\in\O'\}\q\q\q\q\q\q d(\O',\O''):=\inf\{d(x,y):\,x\in\O',\,y\in\O''\}.$$
We will denote as $B_r(x)$ the open ball of radius $r$ centered at $p$. The symbol $\fint_{\O'}f:=\fr{1}{|\O'|}\int_{\O'}f$ will stand for the average of $f\in L^1(\O')$ on some $\O'\sub\O$.\\
The letter $C$ will denote large constant which may vary among different formulas and lines.\\

\section{Compactness issues}\

This section is devoted to the proof of the following concentration-compactness result with quantization of blow-up limits:

\bprop\label{comp}${}$\\
Let $(u_n)_{n\in\N}$ be a sequence of solutions to \eqref{ks} with $\rho_n\us{n\to+\infty}\to\rho$ and $-\b\ne\la_i$ for all $i$'s.\\
Then, up to sub-sequences, one of the following alternatives occur:
\begin{itemize}
\item[(Compactness)] $(u_n)_{n\in\N}$ is compact in $\ol H^1(\O)$;
\item[(Concentration)] The blow-up set $\mcal S$, defined by
$$\mcal S:=\l\{x\in\ol\O:\,\ex\,x_n\us{n\to+\infty}\to x\tx{ such that }u_n(x_n)\us{n\to+\infty}\to+\infty\r\},$$
is non-empty and finite.\\
Moreover, $\rho_n\fr{e^{u_n}}{\int_\O e^{u_n}}\us{n\to+\infty}\wk\sum_{x\in\mcal S}\s(x)\d_x$ as measures, with
\bequ\label{sigma}
\s(x):=\lim_{r\to0}\lim_{n\to+\infty}\rho_n\fr{\int_{B_r(x)}e^{u_n}}{\int_\O e^{u_n}}=\l\{\barr{ll}4\pi&\tx{if }x\in\mcal S\cap\pa\O\\8\pi&\tx{if }x\in\mcal S\cap\O\earr\r..
\eequ
\end{itemize}
In particular, if $\rho\nin4\pi\N$ then (Compactness) occurs.
\eprop\

When \emph{Compactness} occurs, a standard consequence is the following: since the set of solutions is compact, then its energy is uniformly bounded from above, hence the whole space $\ol H^1(\O)$ can be retracted on a suitable sublevel containing all solutions.

\bcor\label{subl}${}$\\
Assume $\rho\nin4\pi\N$ and $-\b\ne\la_i$ for all $i$'s.\\
Then, there exists $L>0$ such that $\mcal J_{\b,\rho}^L$ is a deformation retract of $\ol H^1(\O)$. In particular, it is contractible.
\ecor\

\bpf${}$\\
Since we can write the energy functional as $\mcal J_{\b,\rho}(u)=\fr{\|u\|_{\ol H^1(\O)}^2}2-K_1(u)-\la K_2(u)$, we are in position to apply Proposition $1.1$ in \cite{luc}. From Proposition \ref{comp}, there are no solutions $u_n$ to \eqref{ks} with $L\le\mcal J_{\b,\rho_n}\le L+1$, if $L$ is large enough; therefore, arguing as in \cite{luc}, $\mcal J_{\b,\rho}^L$ is a deformation retract of $\mcal J_{\b,\rho}^{L+1}$. Being $L$ arbitrary, we find that $\mcal J_{\b,\rho}^L$ is a retract of the whole space $\ol H^1(\O)$.
\epf\

Proposition \ref{comp} is rather classical for Liouville-type equations like \eqref{ks}. It was first given by \cite{bremer} and, in the case of Neumann conditions, in \cite{ww}. With respect to the latter reference, the presence of the extra term $\b u$ in the linear part, which may cause the maximum principle to fail, can be dealt by just moving it to the right-hand side.\\
A key tool is a \qm{minimal mass} lemma. The proof in \cite{ww} also works in the case $f_n\not\eq0$.

\blem(\cite{ww}, Lemma $3.2$)\label{minmas}${}$\\
Let $(v_n)_{n\in\N}$ be a sequence of solutions to 
\bequ\label{p0}
\l\{\bll-\D u_n=\rho_n\l(\fr{e^{u_n}}{\int_\O e^{u_n}}-\fr{1}{|\O|}\r)+f_n&\tx{in }\O\\\pa_\nu u_n=0&\tx{on }\pa\O\\\earr\r.
\eequ
with $(f_n)_{n\in\N}$ bounded in $L^q(\O)$ for some $q>1$ with $\int_\O f_n=0$ $\fa n\in\N$.\\
Then, there exists $\s_0=\s_0(\O)>0$ such that if $\limsup_{r\to0}\rho_n\fr{\int_{B_r(x)}e^{v_n}}{\int_\O e^{v_n}}\le\s_0$ for all $x\in\ol\O$, then $(u_n)_{n\in\N}$ is compact in $\ol H^1(\O)$.
\elem\

Roughly speaking, the idea to prove Proposition \ref{comp} will be the following. If \emph{Concentration} occurs, then we have blow-up at a finite number of points, thanks to Lemma \ref{minmas}. Then, the local mass at each blow-up point $x\in\mcal S$ is found via a Poho\v zaev identity based on the asymptotic behavior of solution, and this fact excludes the presence of a residual mass.

\bpf[Proof of Proposition \ref{comp}]${}$\\
Let $(u_n)_{n\in\N}$ be a sequence of solutions with $\sup_\O u_n\le C$. Then, Jensen's inequality gives $\int_\O e^{u_n}\ge|\O|e^{\fint_\O u_n}=|\O|$, therefore
$$|-\D u_n+\b u_n|\le\rho_n\l(\l|\fr{e^{u_n}}{\int_\O e^{u_n}}\r|+\fr{1}{|\O|}\r)\le(\rho+1)\l(\fr{e^C}{|\O|}+\fr{1}{|\O|}\r)$$
is uniformly bounded, hence by standard regularity $u_n$ is bounded in $W^{2,2}(\O)$ and compact in $H^1(\O)$.\\
Suppose now $\sup_\O u_n\us{n\to+\infty}\to+\infty$, namely $\mcal S\ne\es$. This time, we just have $\|-\D u_n+\b u_n\|_{L^1(\O)}\le2(\rho+1)$; since $-\b$ is not an eigenvalue of $-\D$, this gives $\|\n u_n\|_{L^q(\O)}+\|u_n\|_{L^q(\O)}\le C$ for any $q<2$.\\
Therefore, $u_n$ will solve \eqref{p0} with $f_n=-\b u_n\in L^q(\O)$, hence we can apply Lemma \ref{minmas} to get
$$|\mcal S|\s_0\le\sum_{x\in\mcal S_i}\s(x)\le\rho.$$
This means that $\mcal S$ is finite and we easily get $\rho_n\fr{e^{u_n}}{\int_\O e^{u_n}}\us{n\to+\infty}\wk\sum_{x\in\mcal S}\s(x)\d_x+f$ for some $f\in L^1(\O)\cap L^\infty_{\mrm{loc}}(\O\sm\mcal S)$, while $u_n\us{n\to+\infty}\to\sum_{x\in\mcal S}\s(x)G_x+w$ in $W^{1,q}(\O)\cap C^{1,\a}_{\mrm{loc}}(\O\sm\mcal S)$ for $q<2$, $\a<1$, with $G_x$ and $w$ solving, respectively,
$$\l\{\bll-\D G_x+\b G_x=\d_x-\fr{1}{|\O|}&\tx{in }\O\\\pa_\nu G_x=0&\tx{on }\pa\O\\\earr\r.\q\q\q\l\{\bll-\D w+\b w=f-\fint_\O f&\tx{in }\O\\\pa_\nu w=0&\tx{on }\pa\O\\\earr\r..$$
We need to show that $f\eq0$, which will be done arguing as in \cite{batman} (Lemmas $2.2$ and $2.3$). If $f\not\equiv0$, then one easily sees that $f=Ve^w$, with $V:=\fr{\rho}{\ub{\lim_{n\to+\infty}\int_\O e^{u_n}}_{\neq+\infty}}e^{\sum_{x\in\mcal S}\s(x)G_x}\in L^\infty_{\mrm{loc}}(\O\sm\mcal S)$ whereas, due to the behavior of $G_x$, $V\sim|\cd-x|^{-\fr{\s(x)}{2\pi}}$ if $x\in\mcal S\cap\O$ and $V\sim|\cd-x|^{-\fr{\s(x)}\pi}$ if $x\in\mcal S\cap\pa\O$.\\
Now, since $-\D w\ge-\b w-\fint_\O f\in L^q(\O)$, then $w$ is bounded from below, namely $V\le Cf\in L^1(\O)$, which in particular means $\s(x)<4\pi$ for any $x\in\mcal S$. This contradicts \eqref{sigma}, which can be deduced arguing as in \cite{ww} (Lemma $3.4$), hence it must be $f\eq0$.\\
Finally, if \emph{Concentration} occurs then
$$\rho=\lim_{n\to+\infty}\int_{\O}\rho_n\fr{e^{u_n}}{\int_\O e^{u_n}}=\lim_{n\to+\infty}\int_{\O}\sum_{x\in\mcal S}\s(x)\d_x=\sum_{x\in\mcal S}\s(x)=\sum_{x\in\mcal S\cap\pa\O}4\pi+\sum_{x\in\mcal S\cap\O}8\pi\in4\pi\N.$$

\epf\

\section{Analysis of sublevels}\

In this Section, which is the largest of the paper, we will study topologically the energy sublevels of $\mcal J_{\b,\rho}$.\\
In the first sub-section, we will introduce a topological space which will be later \qm{compared} to energy sublevels, and we will compute some of its homology groups. Then, we will construct maps from this topological space to low sublevels and vice-versa and we will deduce that $\mcal J_{\b,\rho}^{-L}$ has non-trivial homology.\\

\subsection{The space $(\O_\pa)_{K,I}$ and its homology}\

Let us introduce a set of \emph{barycenters} on $\ol\O$, namely a set of finitely-supported probability measures on $\ol\O$.\\
With respect to most previous works, we will not give a constraint on the cardinality of the support, basically because points in $\O$ and $\pa\O$ have to be treated differently, for reasons which will be discussed in the forthcoming sub-sections. Roughly speaking, points in the interior will count twice as much as points in the boundary.

$$(\O_\pa)_K:=\Cup_{l=0}^{\l\lfloor\fr{K}2\r\rfloor}\O_{l,K-2l}\q\q\q\O_{l,m}:=\l\{\sum_{k=1}^lt_k\d_{x_k}+\sum_{k'=1}^mt'_{k'}\d_{x'_{k'}};\,\sum_{k=1}^lt_k+\sum_{k'=1}^mt'_{k'}=1,\,x_k\in\O,x'_{k'}\in\pa\O\r\}$$\

On such spaces, we will consider the distance induced by the $\mrm{Lip}'$ norm, that is the norm on the space of signed measures induced by duality with Lipschitz functions:

$$\|\mu\|_{\mrm{Lip}'\l(\ol\O\r)}:=\sup_{h\in\mrm{Lip}\l(\ol\O\r),\|h\|_{\mrm{Lip}\l(\ol\O\r)}\le1}\l|\int_\O h\mrm d\mu\r|.$$\

As a first result, we see that such barycenters spaces are Euclidean deformation retracts.\\
The proof has been given in \cite{ndi} in the case when $\O$ is replaced by a $4$-dimensional compact manifold with boundary, but the same proof holds in any dimension and in particular for planar domains.

\blem(\cite{ndi}, Lemma $4.10$)\label{retr}${}$\\
There exist $\epsilon_0>0$ and a continuous retraction
$$\wt\Psi:\l\{\mu\in\mcal M\l(\ol\O\r):\,d_{\mrm{Lip}'\l(\ol\O\r)}(\mu,(\O_\pa)_K)<\epsilon_0\r\}\to(\O_\pa)_K.$$
\elem\

Among all the \qm{layers} which compose $(\O_\pa)_K$, a special role will be played by the first one, consisting of measures supported on $\pa\O$. In \cite{ndi} it was shown that it is a deformation retract within $(\O_\pa)_K$; again, their proof is also valid in our case.

\blem(\cite{ndi}, Proposition $4.5$)\label{defret}${}$\\
For any $K\ge1$ the set $\O_{0,K}=(\pa\O)_K\sub(\O_\pa)_K$ is a deformation retract of some its open neighborhood $U$ in $(\O_\pa)_K$.
\elem\

Since $\pa\O$ is homotopically equivalent to a disjoint union of $g$ circles, we can use a result from \cite{dlr} to compute the homology of $\O_{0,K}$.

\blem(\cite{dlr}, Proposition $5.1$)${}$\\
The homology groups of $\O_{0,K}=(\pa\O)_K$ are given by
$$\wt H_q((\pa\O)_m)=\l\{\bll\Z^{\bin{g+q-K+1}g\bin{g}{2K-q-1}}&\max\{K-1,2K-g-1\}\le q\le2K-1\\0&q<\max\{K-1,2K-g-1\},q>2K-1\earr\r.$$
\elem\

The space $\l(\O_\pa\r)_K$ will be used in the analysis of sublevels to express the fact that, if $\mcal J_{\b,\rho}(u)\ll0$, then $u$ may concentrates at a finite number of points.\\
Anyway, it $\mcal J_{\b,\rho}$ is very low, it may also happen that the projection $\Pi_I$ on the space of negative eigenvalues for $-\D+\b$, defined by \eqref{pi}, is very large in norm. Since $\Pi_Iu\in\R^I$, this naturally leads to consider, after a normalization, the sphere $\S^{I-1}$ to deal with phenomena.\\
As anticipated in the introduction, the alternative between concentration and large $\Pi_I$ will be modeled by the join \eqref{join}, therefore we will be interested in the following space:

\bequ\label{oki}
(\O_\pa)_{K,I}=\l(\O_\pa\r)_K\star\S^{I-1}.
\eequ\

\brem\label{joineq}${}$\\
In \cite{dm}, the authors used a space of the kind $\fr{X\x\B^{I-1}}\sim$, with $\B^{I-1}$ indicating the unit ball in $\R^{I-1}$ and $\sim$ defined by $(x,y)\sim(x',y)$ for any $x,x'\in X$ and $y\in\S^{I-1}$.\\
Actually, this space is homeomorphic to the join $X\star\S^{I-1}$, with one possible homeomorphism given by the map
$$X\star\S^{I-1}\ni(x,y,t)\q\q\q\ltot\q\q\q(x,ty)\in\fr{X\x\B^{I-1}}\sim.$$
\erem\

Since we are interested in the homology of $(\O_\pa)_{K,I}$, we will use a well-known result concerning the homology of a join.

\blem(\cite{hat}, Theorem $3.21$)\label{omjoin}${}$\\
Let $X$ and $Y$ be two CW-complexes and $X\star Y$ their join as defined by \eqref{join}.\\
Then, its homology group are
$$\wt H_q(X\star Y)=\Pl_{q'=0}^q\wt H_{q'}(X)\ox\wt H_{q-q'-1}(Y),$$
where $\wt H_q$ denotes the \emph{reduced homology groups}: $H_q(X)=\l\{\bll\wt H_q(X)\oplus\Z&\tx{if }q=0\\\wt H_q(X)&\tx{if }q\ge1\earr\r.$.\\
\elem\

We are now able to get some information on the homology on $(\O_\pa)_{K,I}$. In particular, we will compute its maximal dimensional homology group, which is not trivial.

\bprop\label{omo}${}$\\
Let $g$ be the genus of $\O$.\\
The homology groups of $(\O_\pa)_{K,I}$ satisfy
$$\wt H_{2K+I-1}((\O_\pa)_{K,I})=\Z^{\bin{K+g}g}.$$
\eprop\

\bpf${}$\\
We start with the case $I=0$. If $K=1$, then $(\O_\pa)_1=\pa\O$ so its homology is computed immediately; otherwise, we write $(\O_\pa)_K=U\cup V$, with $U$ as in Lemma \ref{defret} and $V:=(\O_\pa)_K\sm(\pa\O)_K$. The space $V=\Cup_{l=0}^{\l\lfloor\fr{K}2\r\rfloor}\O_{l,K-2l}$ is a stratified set whose maximal dimension equals $2K-3$, and the same holds true for $U\cap V$, hence $\wt H_q(U\cap V)=H(U\cap V)=0$ for any $q\ge2K-2$. Therefore, the Mayer-Vietoris exact sequence gives:
$$0=\wt H_{2K-1}(U\cap V)\to\wt H_{2K-1}(U)\oplus\wt H_{2K-1}(V)\to\wt H_{2K-1}((\O_\pa)_K)\to\wt H_{2K-2}(U\cap V)=0,$$
that is
$$\wt H_{2K-1}((\O_\pa)_K)=\wt H_{2K-1}(U)\oplus H_{2K-1}(V)=\wt H_{2K-1}((\pa\O)_K)=Z^{\bin{K+g}g}.$$

Finally, if $I\ge1$, then Lemma \ref{omjoin} gives
\beqy
\wt H_{2K+I-1}((\O_\pa)_{K,I})&=&\wt H_{2K+I-1}\l((\O_\pa)_K\star\S^{I-1}\r)\\
&=&\Pl_{q'=0}^{2K+I-1}\wt H_{q'}((\O_\pa)_K)\ox\wt H_{2K+I-q'-2}\l(\S^{I-1}\r)\\
&=&H_{2K-1}((\O_\pa)_K)\\
&=&\Z^{\bin{K+g}g}.
\eeqy

\brem${}$\\
As pointed out by the referee, one can compute the Euler characteristic of $(\O_\pa)_K$ using the results in \cite{kk,akn}:
$$\chi((\O_\pa)_K)=\chi\l(\O_{\l\lfloor\fr{K}2\r\rfloor,0}\r)=1-\fr{1}{\l\lfloor\fr{K}2\r\rfloor!}\prod_{k=1}^{\l\lfloor\fr{K}2\r\rfloor}(k-\chi(\O))=1-\bin{\l\lfloor\fr{K}2\r\rfloor+g-1}{g-1}.$$
\erem\

\epf\

\subsection{The map $\Phi^\La:(\O_\pa)_{K,I}\to\mcal J_{\b,\rho}^{-L}$}\

In this Subsection we will build a map from the space $(\O_\pa)_{K,I}$, whose properties have just been discussed, into a suitably low energy sublevel $J_{\b,\rho}^{-L}$.\\
Precisely, we will construct a \emph{family} $\Phi^\La$ of maps with $\mcal J_{\b,\rho}\l(\Phi^\La\r)\us{\La\to+\infty}\to-\infty$ uniformly, so that for $\La$ large enough the image of $\mcal J_{\b,\rho}$ is contained in $\mcal J_{\b,\rho}^{-L}$. The choice of $L$, hence of $\La$, will be made in the following subsection.\\
Consistently with the previous discussion, the family of test functions defined by $\Phi^\La$ will have the following property: as $\La$ goes to $0$, either it concentrates at a finite number of points, according to the definition of $(\O_\pa)_K$ (if $t\ne1$), or its projection $\Pi_I$ will be large (if $t\ne0$).

\bprop\label{test}${}$\\
Let $(\O_\pa)_{K,I}$ be defined by \eqref{oki} and let $\Phi^\La:(\O_\pa)_{K,I}\to\ol H^1(\O)$ be defined, for $L\gg0$, in the following way:
$$\z=(\mu,\si,t)=\l(\sum_kt_k\d_{x_k},(\si_1,\ds,\si_I),t\r)\q\q\q\longmapsto\q\q\q\Phi^\La(\z):=\phi^{\La(1-t)}-\fint_\O\phi^{\La(1-t)}+\psi^{\La t}$$
\beqy
\phi^{\La(1-t)}=\phi^{\La(1-t)}(\mu)&:=&\log\sum_k\fr{t_k}{(1+(\La(1-t))^2|\cd-x_k|^2)^2}\\
\psi^{\La t}=\psi^{\La t}(\si)&:=&\sqrt{\log^+(\La t)}\sum_{i=1}^I\si_i\ph_i
\eeqy
If $\rho>4K\pi$ and $\b<-\la_I$, then $\mcal J_{\b,\rho}\l(\Phi^\La(\z)\r)\us{\La\to+\infty}\to-\infty$ independently on $\z$.\\
\eprop\

To prove this result, we will estimate separately the three parts defining $\mcal J_{\b,\rho}$: the Dirichlet integral, the $L^2$ norm and the nonlinear term. Each estimate is contained in a separate lemma.

\blem\label{grad}${}$\\
Let $\Phi^\La$ be as in Proposition \ref{test}.\\
Then,
$$\int_\O\l|\n\Phi^\La(\z)\r|^2\le16K\pi\log^+(\La(1-t))+\la_I\log^+(\La t)+C\sqrt{\log^+(\La t)}.$$
\elem\

\bpf${}$\\
Since, by definition, we have $$\l|\n\Phi^\La(\z)\r|^2=\l|\n\phi^{\La(1-t)}\r|^2+2\n\phi^{\La(1-t)}\cd\n\psi^{\La t}+\l|\n\psi^{\La t}\r|^2,$$
we will suffice to show the following estimates:
\beqa
\label{phi2}\int_\O\l|\n\phi^{\La(1-t)}\r|^2&\le&16K\pi\log^+(\La(1-t))+C;\\
\label{phipsi}\int_\O\n\phi^{\La(1-t)}\cd\n\psi^{\La t}&\le&C\sqrt{\log^+(\La t)};\\
\label{psi2}\int_\O\l|\n\psi^{\La t}\r|^2&\le&\la_I\log^+(\La t).
\eeqa

The first estimate can be obtained similarly as \cite{mal1} (Proposition $4.2$), the main difference being that we have to take care of the points $x_k$ lying on $\pa\O$. The estimate is trivial if $\La(1-t)$ is bounded from above, otherwise we get:
\beqy
\l|\n\phi^{\La(1-t)}\r|&=&\l|\fr{\sum_k\fr{-4t_k(\La(1-t))^2(\cd-x_k)}{(1+(\La(1-t))^2|\cd-x_k|^2)^3}}{\sum_k\fr{t_k}{\l(1+(\La(1-t))^2|\cd-x_k|^2\r)^2}}\r|\\
&\le&\fr{\sum_k\fr{4t_k(\La(1-t))^2|\cd-x_k|^2}{(1+(\La(1-t))^2|\cd-x_k|^2)^3}}{\sum_k\fr{t_k}{\l(1+(\La(1-t))^2|\cd-x_k|^2\r)}}\\
&\le&\max_k\fr{4(\La(1-t))^2|\cd-x_k|}{1+(\La(1-t))^2|\cd-x_k|^2}\\
&\le&\min\l\{4\La(1-t),\fr{4}{\min_k|\cd-x_k|}\r\}.
\eeqy
Now, we divide $\O$ in some regions $\O_k$ depending on which point $x_k$ is the closest:
$$\O_k:=\l\{x\in\O:\,|x-x_k|=\min_{k'}|x-x_{k'}|\r\};$$
therefore, we get
\beqy
\int_\O\l|\n\phi^{\La(1-t)}\r|^2&\le&\sum_k\int_{\O_k}\l|\n\phi^{\La(1-t)}\r|^2\\
&\le&\sum_k\int_{\O_k\sm B_\fr{1}{\La(1-t)}}\fr{16}{|\cd-x_k|^2}+\sum_k\int_{B_\fr{1}{\La(1-t)}}16(\La(1-t))^2\\
&\le&16\sum_k\int_{\O_k\sm B_\fr{1}{\La(1-t)}}\fr{1}{|\cd-x_k|^2}+C.
\eeqy
To evaluate the last integral we distinguish the cases $x_k\in\O$ and $x_k\in\pa\O$. In the former, we are basically integrating the function $\fr{1}{|\cd|^2}$ on an annulus whose internal radius $\fr{1}{\La(1-t)}$ is shrinking, plus negligible terms, hence its asymptotical value will be $2\pi\log(\La(1-t))$. On the other hand, if $x_k\in\pa\O$, then we are actually integrating on a domain asymptotically resembling a half-annulus with its internal radius shrinking, therefore we only get half of before, namely $\pi\log(\La(1-t))$.\\
Following these considerations, we get \eqref{phi2}:
\beqy
\l|\n\phi^{\La(1-t)}\r|&\le&16\l(\sum_{x_k\in\O}(2\pi\log(\La(1-t))+C)+\sum_{x_k\in\pa\O}\pi(\log(\La(1-t))+C)\r)+C\\
&\le&16K\log(\La(1-t))+C.
\eeqy\

Concerning \eqref{phipsi}, by the construction of $\psi^\La$ we get
\beqy
\int_\O\n\phi^{\La(1-t)}\cd\n\psi^{\La t}&=&\sqrt{\log^+(\La t)}\sum_{i=1}^Is_i\int_\O\n\phi^{\La(1-t)}\cd\n\ph_i\\
&=&\sqrt{\log^+(\La t)}\sum_{i=1}^Is_i\la_i\int_\O\l(\phi^{\La(1-t)}-\fint_\O\phi^{\La(1-t)}\r)\ph_i\\
&\le&C\sqrt{\log^+(\La t)}\sum_{i=1}^Is_i\la_i\sqrt{\int_\O\l(\phi^{\La(1-t)}-\fint_\O\phi^{\La(1-t)}\r)^2}\sqrt{\int_\O\ph_i^2}\\
&\le&C\sqrt{\log^+(\La t)}\sqrt{\int_\O\l(\phi^{\La(1-t)}-\fint_\O\phi^{\La(1-t)}\r)^2},
\eeqy
therefore we suffice to show that the last integral is uniformly bounded.\\
To this purpose, we first estimate the average of $\phi^{\La(1-t)}$:
\beqa
\nonumber\fint_\O\phi^{\La(1-t)}&=&\fint_\O\log\fr{1}{(1+(\La(1-t))^2\min_k|\cd-x_k|^2)^2}+O(1)\\
\nonumber&=&\fint_{\O\sm\Cup_kB_\fr{1}{\La(1-t)}}\log\fr{1}{(\La(1-t)\min_k|\cd-x_k|)^4}+O(1)\\
\label{av}&=&-4\log^+(\La(1-t))+O(1).
\eeqa
Now, \eqref{phipsi} will follow by estimating $\phi^{\La(1-t)}+4\log^+(\La(1-t))$ in $L^2(\O)$, which can be done similarly as before:
\beqa
\nonumber\int_\O\l(\phi^{\La(1-t)}+4\log^+(\La(1-t))\r)^2&=&\int_\O\l(\log\sum_k\fr{t_k\max\{1,\La(1-t)\}}{(1+(\La(1-t))^2|\cd-x_k|^2)^2}\r)^2\\
\nonumber&&\int_\O\log\l(\sum_k\fr{t_k}{|\cd-x_k|^4}\r)^2\\
\nonumber&\le&\int_\O\log\fr{1}{\min_k|\cd-x_k|^8}\\
\label{l2}&\le&C.
\eeqa\

Finally, \eqref{psi2} follows easily by the properties of the $\phi_i$'s:
\beqy
\int_\O\l|\n\psi^{\La t}\r|^2&=&\log^+(\La t)\int_\O\l|\sum_{i=1}^I\si_i\n\ph_i\r|^2\\
&=&\log^+(\La t)\sum_{i=1}^I\si_i^2\int_\O|\n\ph_i|^2\\
&=&\log^+(\La t)\sum_{i=1}^I\si_i^2\la_i\\
&\le&\log^+(\La t)\la_I.\\
\eeqy
\epf\

\blem\label{sq}${}$\\
Let $\Phi^\La$ be as in Proposition \ref{test}.\\
Then,
$$\int_\O\Phi^\La(\z)^2\le\log^+(\La t)+C\sqrt{\log^+(\La t)}.$$
\elem\

\bpf${}$\\
By expanding the square of the sum and using \eqref{av}, \eqref{l2}, we have
\beqy
\int_\O{\Phi^{\La}(\z)}^2&=&\int_\O\l(\phi^{\La(1-t)}-\fint_\O\phi^{\La(1-t)}\r)^2+2\int_\O\l(\phi^{\La(1-t)}-\fint_\O\phi^{\La(1-t)}\r)\psi^{\La t}+\int_\O\l(\psi^{\La t}\r)^2\\
&\le&\int_\O\l(\phi^{\La(1-t)}-\fint_\O\phi^{\La(1-t)}\r)^2+2\sqrt{\int_\O\l(\phi^{\La(1-t)}-\fint_\O\phi^{\La(1-t)}\r)^2}\sqrt{\int_\O\l(\psi^{\La t}\r)^2}+\int_\O\l(\psi^{\La t}\r)^2\\
&\le&C\l(1+\sqrt{\int_\O\l(\psi^{\La t}\r)^2}\r)+\int_\O\l(\psi^{\La t}\r)^2.
\eeqy
Therefore, we only need a suitable estimate for $\psi^{\La t}$, which in turn follows from the very definition of the $\ph_i$'s:
$$\int_\O\l(\psi^{\La t}\r)^2=\log^+(\La t)\int_\O\l(\sum_{i=1}^I\si_i\ph_i\r)^2=\log^+(\La t)\sum_{i=1}^I\si_i^2\int_\O\ph_i^2=\log^+(\La t)\sum_{i=1}^I\si_i^2\le\log^+(\La t).$$
By combining the two estimates the proof is complete. 
\epf\

\blem\label{exp}${}$\\
Let $\Phi^\La$ be as in Proposition \ref{test}.\\
Then,
$$\log\int_\O e^{\Phi^\La(\z)}\ge2\log^+(\La(1-t))-C\sqrt{\log^+(\La t)}.$$
\elem\

\bpf${}$\\
We first notice that, since $\psi^{\La t}$ belongs to a finite-dimensional space, all of its norms are equivalent, and in particular $L^2$ and $L^\infty$, therefore by Lemma \eqref{sq},
$$\l|\psi^{\La t}\r|\le\l\|\psi^{\La t}\r\|_{L^\infty(\O)}\le C\l\|\psi^{\La t}\r\|_{L^2(\O)}\le C\sqrt{\log^+(\La t)}.$$
Moreover, in view of the asymptotical behavior \eqref{av} of the average of $\phi_{\la(1-t)}$, we are reduce to show that $\log\int_\O e^{\phi^{\La(1-t)}}\ge-2\log^+(\La(1-t))-C$; this follows from the simple calculations:
\beqy
\int_\O e^{\phi^{\La(1-t)}}&=&\sum_kt_k\int_\O\fr{1}{(1+(\La(1-t))^2|\cd-x_k|^2)^2}\\
&\ge&\sum_kt_k\int_{B_\fr{1}{\La(1-t)}(x_k)}\fr{1}{(1+(\La(1-t))^2|\cd-x_k|^2)^2}\\
&\ge&\sum_kt_k\int_{B_\fr{1}{\La(1-t)}(x_k)}\fr{1}2\\
&\ge&\fr{C}{\max\{1,\La(1-t)\}^2}.
\eeqy
\epf\

By putting together these three lemmas, Proposition \eqref{test} may be proved easily.

\bpf[Proof of Proposition \ref{test}]${}$\\
By Lemmas \ref{grad}, \ref{sq}, \ref{exp}, we get
\beqy
\mcal J_{\b,\rho}\l(\Phi^\La(\z)\r)&=&\fr{1}2\int_\O\l|\n\Phi^\La(\z)\r|^2+\fr{\b}2\int_\O\Phi^\La(\z)^2-\rho\log\int_\O e^{\Phi^\La(\z)}\\
&\le&(8K\pi-2\rho)\log^+(\La(1-t))+\fr{\la_I+\b}2\log^+(\La t)+C\sqrt{\log^+(\La t)}\\
&\le&-\min\l\{2\rho-8K\pi,-\fr{\la_I+\b}2\r\}\max\l\{\log^+(\La(1-t)),\log^+(\La t)\r\}+C\sqrt{\log^+(\La t)}\\
&\le&-\min\l\{2\rho-8K\pi,-\fr{\la_I+\b}2\r\}\log\fr{\La}2+C\sqrt{\log\La}\\
&\us{\La\to+\infty}\to&-\infty,
\eeqy
uniformly on $\z\in(\O_\pa)_{K,I}$.
\epf\

\subsection{The map $\Psi:\mcal J_{\b,\rho}^{-L}\to(\O_\pa)_{K,I}$}\

We will now show the existence of \qm{counterpart} to the map $\Phi$ defined in the previous subsection.\\
Precisely, we will build a map $\Psi$ from a low sub-level $\mcal J_{\b,\rho}^{-L}$ to $(\O_\pa)_{K,I}$ which is somehow compatible with $\Phi$, in the sense that their composition is homotopically equivalent to the identity on $(\O_\pa)_{K,I}$.

\bprop\label{psi}${}$\\
Let $(\O_\pa)_{K,I}$ be defined by \eqref{oki}.\\
If $\rho<4(K+1)\pi$ and $\b>-\la_{I+1}$, then there exists $L\gg0$ and a map $\Psi:\mcal J_{\b,\rho}^{-L}\to(\O_\pa)_{K,I}$.\\
Moreover, if $\rho>4K\pi$ and $\b<-\la_I$, then there exists a map $\Phi:(\O_\pa)_{K,I}\to\mcal J_{\b,\rho}^{-L}$ such that the composition $\Psi\c\Phi$ is homotopically equivalent to the identity on $(\O_\pa)_{K,I}$.
\eprop\

The existence of such maps $\Phi$ and $\Psi$ easily gives, via the functorial properties of homology, the following information on the homology groups of energy sublevels.

\bcor\label{omosub}${}$\\
Under the assumptions of Propositions \ref{test} and \ref{psi}, the homology groups of the sublevel $\mcal J_{\b,\rho}^{-L}$ satisfy
$$\Z^{\bin{K+g}g}\inc\wt H_{2K+I-1}\l(\mcal J_{\b,\rho}^{-L}\r).$$
\ecor\

The main tool to prove Proposition \ref{psi} is a so-called \emph{improved} Moser-Trudinger inequality. Roughly speaking, such inequalities state that, under some \emph{spreading} conditions on $u$, the best constant in the classical Moser-Trudinger inequality can be improved.\\
We recall here the well-known Moser-Trudinger inequalities, in two forms depending whether we consider only compactly supported function or also function which may touch the boundary. We stress that, in the two cases, the constant multiplying the Dirichlet integral is different.

\bprop(\cite{mos}, Theorem $2$, \cite{cy88}, Corollary $2.5$)${}$\\
There exists $C>0$ such that for any $u\in\ol H^1(\O)$
\bequ\label{mt1}
\log\int_\O e^u\le\fr{1}{8\pi}\int_\O|\n u|^2+C.
\eequ
If instead $u\in H^1_0(\O)$, then
\bequ\label{mt2}
\log\int_\O e^u\le\fr{1}{16\pi}\int_\O|\n u|^2+C.
\eequ
\eprop\

We will now prove the improved Moser-Trudinger inequality, a classical result in variational Liouville-type problems (see \cite{dm,mal1,cm,bjmr,bat1,bat2}).\\
Basically, if $u$ is somehow spread in some regions, then the constant $8\pi$ in \eqref{mt1} can almost be multiplied by an integer number. This time, the integer will depend not only on the number of regions but also on how many of them touch the boundary; moreover, we need to take account of the negative projection $\Pi_I$.\\
To prove such a result, we will take cutoff functions on the regions where $u$ is spread and apply to each cutoff either \eqref{mt1} or \eqref{mt2}. We will also use some splitting in Fourier modes, which we need to deal with $\Pi_I$.

\blem\label{mtimpr}${}$\\
Let $\d>0$, $\{\O_{1k}\}_{k=1}^l$, $\{\O_{2k'}\}_{k'=1}^m$ and $u\in\ol H^1(\O)$ satisfying
$$\barr{lll}
d(\O_{ik},\O_{i'k'})\ge2\d\q\fa\,(i,k)\ne(i',k');&\q\q&d(\O_{1k},\pa\O)\ge\d\q\fa\,k=1,\ds,l;\\
\fr{\int_{\O_{ik}}e^u}{\int_\O e^u}\ge\d\q\fa\,i,k;&\q\q&\|\Pi_Iu\|\le1.
\earr$$
Then, for any $\e>0$ there exists $C=C(\e,\d,\b,I,l,m)>0$ such that
$$\log\int_\O e^u\le\fr{1+\e}{8\pi(2l+m)}\int_\O\l(|\n u|^2+\b u^2\r)+C.$$
\elem\

\bpf${}$\\
First of all, for any $i,k$ we take cutoff functions $\eta_{ik}\in\mrm{Lip}\l(\ol\O\r)$ satisfying
$$0\le\eta_{ik}\le1\q\tx{in }\ol\O_{ik}\q\q\q\q\eta_{ik}|_{\O_{ik}}\eq1\q\q\q\q\mrm{spt}(\eta_{ik})\sub B_\d(\O_{ik})\q\q\q\q|\n\eta_{ik}|\le\fr{1}\d\q\tx{in }B_\d(\O_{ik})\sm\O_{ik}.$$
Then, we split $u=u_1+u_2+u_3$ via truncation in Fourier modes:
$$u_1=\sum_{i=1}^I\l(\int_\O u\ph_i\r)\ph_i\q\q\q\q\q\q u_2=\sum_{i=I+1}^{N_\e-1}\l(\int_\O u\ph_i\r)\ph_i\q\q\q\q\q\q u_3=\sum_{i=N_\e}^{+\infty}\l(\int_\O u\ph_i\r)\ph_i,$$
with $N_\e$ so large that
\bequ\label{la}
\fr{1}{\la_{N_\e}}\l(1+\fr{1}\e\r)\fr{1}{\d^2}\le\e\q\q\q\q\q\q\la_{N_\e}\le(1+\e)(\la_{N_\e}+\b).
\eequ
By applying \eqref{mt1} to each $\eta_{1k}u-\ol{\eta_{1k}u}$ we get:
\beqy
\log\int_\O e^u&\le&\log\int_{\O_{1k}}e^u+\log\fr{1}\d\\
&\le&\log\int_\O e^{\eta_{1k}u}+\log\fr{1}\d\\
&\le&\|\eta_{1k}u_1\|_{L^\infty(\O)}+\|\eta_{1k}u_2\|_{L^\infty(\O)}+\log\int_\O e^{\eta_{1k}u_3}+\log\fr{1}\d\\
&\le&\|u_1\|_{L^\infty(\O)}+\|u_2\|_{L^\infty(\O)}+\ol{\eta_{1k}u_3}+\fr{1}{8\pi}\int_\O|\n(\eta_{2k}u_3)|^2+C.
\eeqy
Similarly, since $\eta_{2k}u\in H^1_0(\O)$, we can apply \eqref{mt2}:
\beqy
\log\int_\O e^u&\le&\log\int_{\O_{2k}}e^u+\log\fr{1}\d\\
&\le&\log\int_\O e^{\eta_{2k}u}+\log\fr{1}\d\\
&\le&\|\eta_{2k}u_1\|_{L^\infty(\O)}+\|\eta_{2k}u_2\|_{L^\infty(\O)}+\log\int_\O e^{\eta_{2k}u_3}+\log\fr{1}\d\\
&\le&\|u_1\|_{L^\infty(\O)}+\|u_2\|_{L^\infty(\O)}+\fr{1}{16\pi}\int_\O|\n(\eta_{1k}u_3)|^2+C.
\eeqy
The terms involving $u_1$ and $u_2$ can be estimated because, since each belongs to a finite-dimensional space, all of its norms are equivalent on the respective space. We can use the $L^2$ norm for $u_1$, which is uniformly bounded by hypotheses:
$$\|u_1\|_{L^\infty(\O)}\le C\|u_1\|_{L^2(\O)}=C\|\Pi_Iu\|\le C.$$
As for $u_2$, since we got rid of low Fourier coefficients, we can choose as a norm $\sqrt{\int_\O\l(|\n u_2|^2+\b u_2^2\r)}$; since we took an orthogonal decomposition, we get
\beqy
\|u_2\|_{L^\infty(\O)}&\le&C\sqrt{\int_\O\l(|\n u_2|^2+\b u_2^2\r)}\\
&\le&\e\int_\O\l(|\n u_2|^2+\b u_2^2\r)+C\\
&\le&\e\int_\O\l(|\n u_2|^2+\b u_2^2\r)+\e\int_\O\l(|\n u_3|^2+\b u_3^2\r)+\e\int_\O|\n u_1|^2+C\\
&=&\e\int_\O\l(|\n u|^2+\b u^2\r)-\e\b\|\Pi_Iu\|^2+C\\
&\le&\e\int_\O\l(|\n u|^2+\b u^2\r)+C.
\eeqy
We then estimate the average of $\eta_{1k}u_3$ via Poincaré-Wirtinger inequality:
$$\l|\ol{\eta_{1k}u_3}\r|\le\|u_3\|_{L^1(\O)}\le\|\n u_3\|_{L^2(\O)}\le\e\int_\O|\n u_3|^2+C.$$
Concerning the last term, we expand the square and use the properties of the $\eta_{ik}$'s:
\beqy
\int_\O|\n(\eta_{ik}u_3)|^2&=&\int_\O|\eta_{ik}\n u_3+u_3\n\eta_{ik}|^2\\
&\le&(1+\e)\int_\O\eta_{ik}^2|\n u_3|^2+\l(1+\fr{1}\e\r)\int_\O u_3^2|\n\eta_{ik}|^2\\
&\le&(1+\e)\int_{B_\d(\O_{ik})}|\n u_3|^2+\l(1+\fr{1}\e\r)\fr{1}{\d^2}\int_{B_\d(\O_{ik})} u_3^2.
\eeqy
Since by hypothesis $B_\d(x_{ik})\cap B_\d(x_{i'k'})=\es$ for $(i,k)\ne(i',k')$, then putting together all these estimates and summing on $i=1,2$ and all $k$'s we get
$$(2l+m)\int_\O e^u\le(2l+m)\e\int_\O\l(|\n u|^2+\b u^2\r)+m\e\int_\O|\n u_3|^2+\fr{1+\e}{8\pi}\int_\O|\n u_3|^2+\fr{1}{8\pi}\l(1+\fr{1}\e\r)\fr{1}{\d^2}\int_\O u_3^2+C.$$
At this point, we need the conditions \eqref{la} defining $u_3$: the former gives
$$\l(1+\fr{1}\e\r)\fr{1}{\d^2}\int_\O u_3^2\le\fr{1}{\la_{N_\e}}\l(1+\fr{1}\e\r)\fr{1}{\d^2}\int_\O|\n u_3|^2\le\e\int_\O|\n u_3|^2.$$
on the other hand, the latter implies
\beqy
\int_\O|\n u_3|^2&=&\sum_{j=N_\e}^{+\infty}\la_i\l(\int_\O u\ph_i\r)^2\\
&\le&(1+\e)\sum_{j=N_\e}^{+\infty}(\la_i+\b)\l(\int_\O u\ph_i\r)^2\\
&=&(1+\e)\int_\O\l(|\n u_3|^2+\b u_3^2\r)\\
&\le&(1+\e)\int_\O\l(|\n u_3|^2+\b u_3^2\r)+(1+\e)\int_\O|\n u_1|^2+(1+\e)\int_\O\l(|\n u_2|^2+\b u_2^2\r)\\
&\le&(1+\e)\int_\O\l(|\n u|^2+\b u^2\r)+C.\\
\eeqy
Therefore, we get
\beqy
(2l+m)\int_\O e^u&\le&(2l+m)\e\int_\O\l(|\n u|^2+\b u^2\r)+\l(\fr{1+\e}{8\pi}+\fr{\e}{8\pi}+m\e\r)\int_\O|\n u_3|^2+C\\
&\le&\l((2l+m)\e+(1+\e)\l(\fr{1+\e}{8\pi}+\fr{\e}{8\pi}+m\e\r)\r)\int_\O\l(|\n u|^2+\b u^2\r)+C,
\eeqy
which, up to re-labeling $\e$, concludes the proof.
\epf\

We need a few more technical steps to prove Proposition \ref{psi}.\\
First of all, we need a covering lemma basically saying that, if concentration does not occur, then one has spreading in the sense of Lemma \ref{mtimpr}.

\blem\label{cover}${}$\\
Let $f\in L^1(\O)$ be non-negative a.e., satisfying $\int_\O f=1$ and such that, for any $\e>0$ and $x_{11,},\ds,x_{1l}\in\O,\,x_{21},\ds,x_{2m}\in\pa\O$ with $2l+m\le K$ for some $K\in\N$,
$$\int_{\Cup_{i,k}B_\e(x_{ik})}f<1-\e.$$
Then, there exist $\wt\e=\wt\e(\e,\O),\wt r=\wt r(\e,\O)>0$ and $\wt x_{11},\ds,\wt x_{1\wt l},\wt x_{21},\ds,\wt x_{2\wt m}$ satisfying
$$\barr{lll}
2\wt l+\wt m\ge K+1;&\q\q&d(\wt x_{1k},\pa\O)\ge\wt r\q\fa\,k=1,\ds,\wt l\\
|\wt x_{ik}-\wt x_{i'k'}|\ge4\wt r\q\fa\,(i,k)\ne(i',k')&\q\q&\int_{B_{\wt r}(\wt x_{ik})}f\ge\wt\e\q\fa\,i,k
\earr$$
\elem\

\bpf${}$\\
We will mostly argue as in \cite{dm} (Lemma $2.3$) and \cite{mal1} (Lemma $3.3$), with minor modifications.\\
Fix $\wt r:=\fr{\e}6$ and take the finite cover of $\ol\O$ given by $\{B_{\wt\e}(y_n)\}_{n=1}^N$ for some $L=L_{\wt r,\ol\O}$, then set $\wt\e:=\fr{\e}L$. One easily sees that there exists some $n$ such that $\int_{B_{\wt r}(y_n)}f\ge\wt\e$; up to re-labeling, we can assume that this hold true if and only if $n\le N'$, for some $N'\le N$.\\
Now, we choose recursively the points $\{\wt y_j\}\sub\{y_n\}_{n=1}^{N'}$: we set $\wt y_1:=y_1$ and
$$\O_1:=\l\{\Cup_{n=1}^{N'}B_{\wt r}(y_n):\,|y_n-\wt y_1|<4\wt r\r\}\sub B_{5\wt r}(\wt y_1).$$
If there is some $l_0$ such that $|y_{n_0}-\wt y_1|\ge4\wt r$, then we set $\wt y_2=y_{n_0}$ and
$$\O_2:=\l\{\Cup_{n=2}^{N'}B_{\wt r}(y_n):\,|y_n-\wt y_2|<4\wt r\r\}\sub B_{5\wt r}(\wt y_2).$$
Inductively, we find a finite number of points $\wt y_j$ and closed set $\O_j$.\\
Among the $\wt y_j$'s, some of them will be at a distance less than $\d$ from $\pa\O$; we denote the number of such points as $\wt m$ and the number of the other $\wt y_j$'s as $\wt l$, then we denote the former points as $\wt x_{1k}$ and the latter as $\wt x_{2k'}$, so that $\{\wt y_j\}_j=\l\{\wt x_{11},\ds,\wt x_{1\wt l},\wt x_{21},\ds,\wt x_{2\wt m}\r\}$ and we call $\O_{ik}$ the set $\O_j$ corresponding to the point $\wt x_{ik}$. To complete the proof, we only need to show that $2\wt l+\wt m>K$, since we already verified that the other required properties are satisfied.\\
Assume by contradiction that $2\wt l+\wt m\le K$, set $x_{1k}:=\wt x_{1k}$ for $k=1,\ds,\wt l$ and take, for $k=1,\ds,\wt m$, some $x_{2k'}\in\pa\O$ such that $d\l(x_{2k'},\wt x_{2k}\r)\le\wt r$. Then, by hypothesis,
$$\int_{\O\sm\Cup_{i,k}B_\e(x_{ik})}f\ge\e.$$
However, due to our construction,
$$\Cup_{n=1}^{N'}B_{\wt r}(y_n)\sub\Cup_{i,k}\O_{ik}\sub\Cup_{i,k}B_{5\wt r}\l(\wt x_{ik}\r)\sub\Cup_{i,k}B_\e(x_{ik}),$$
which leads to a contradiction:
$$\int_{\O\sm\Cup_{i,k}B_\e(x_{ik})}f\le\int_{\O\sm\Cup_{n=1}^{N'}B_{\wt r}(y_n)}f\le\int_{\Cup_{n=1}^{N'}B_{\wt r}(y_n)}f\le(N-N')\wt\e<\e$$
\epf\

Now we see that either concentration at a finite number of points or large $\|\Pi_I\|$ must occur in very low sublevels. In fact, if it does not, then by the previous lemma one has spreading and small $\Pi_I$, hence by the improved Moser-Trudinger inequality the energy $\mcal J_{\b,\rho}$ cannot be too low.\\
This is explained in details by the following lemma.

\blem\label{epsl}${}$\\
For any $\e>0$ there exists $L=L(\e)>0$ such that, if $J_{\b,\rho}(u)\le-L$ and $\|\Pi_Iu\|>1$, then there exists $x_{11},\ds,x_{1l}\in\O,\,x_{21},\ds,x_{2m}\in\pa\O$ with $2l+m\le K$ such that
$$\fr{\int_{\Cup_{i,k}B_\e(x_{ik})}e^u}{\int_\O e^u}\ge1-\e.$$
\elem\

\bpf${}$\\
Assume that the statement is not true. Then, there exists $\e>0$ and $(u_n)_{n\in\N}$ such that $\|\Pi_Iu_n\|\le1$, $J_{\b,\rho}(u_n)\us{n\to+\infty}\to-\infty$ and
$$\fr{\int_{\Cup_{i,k}B_\e(x_{ik})}e^{u_n}}{\int_\O e^{u_n}}<1-\e.$$
for any $x_{11},\ds,x_{1l}\in\O,\,x_{21},\ds,x_{2m}\in\pa\O$ satisfying $2l+m\le K$.\\
We apply lemma \ref{cover} to $f=\fr{e^{u_n}}{\int_\O e^{u_n}}$ and we find $\wt\e,\wt r$, not depending on $n$, and $\wt x_{11},\ds,\wt x_{1\wt l},\wt x_{21},\ds,\wt x_{2\wt m}$ as in the lemma.\\
One can easily see that Lemma \ref{mtimpr} can be applied to $\O_{ik}=B_{\wt r}(\wt x_{ik})$ with $\d=\min\{\wt\e,\wt r\}$ and $\e'=\fr{4\pi}\rho\l(2\wt l+\wt m\r)-1$. This leads to the following contradiction:
$$-\infty\us{n\to+\infty}\ot J_{\b,\rho}(u_n)=\fr{4\pi\l(2\wt l+\wt m\r)}{1+\e'}\l(\fr{1+\e'}{8\pi\l(2\wt l+\wt m\r)}\int_\O\l(|\n u_n|^2+\b u_n^2\r)-\rho\log\int_\O e^{u_n}\r)\ge-C$$
\epf\

Since $\fr{e^u}{\int_\O e^u}$ tends to concentrates in very low sublevels, provided $\Pi_I$ is not too large, then it will be very close to an element of $(\O_\pa)_K$. This will be essential to later use the retraction $\wt\Psi$ defined in Lemma \ref{retr}.

\blem\label{deps}${}$\\
For any $\e>0$ there exists $L=L(\e)>0$ such that any $u\in\mcal J_{\b,\rho}^{-L}$ satisfies either of the following condition:
$$d_{\mrm{Lip}'\l(\ol\O\r)}\l(\fr{e^u}{\int_\O e^u},(\O_\pa)_K\r)\le\e\q\q\q\q\q\q\tx{or }\q\q\q\q\q\q\|\Pi_Iu\|>1.$$
\elem\

\bpf${}$\\
Fix $\e>0$, apply Lemma \ref{epsl} with $\fr{\e}3$ and take $L=L\l(\fr{\e}3\r)$ as in the lemma. For any $u\in\mcal J_{\b,\rho}^{-L}$ satisfying $\|\Pi_Iu\|\le1$, take $\mu(u)=\sum_{ik}t_{ik}\d_{x_{ik}}$ with $x_{ik}$ as in the lemma and $t_{ik}$ defined by
$$t_k:=\int_{B_\fr{\e}3(x_{ik})\sm\Cup_{k'\le k-1\tx{ or }i'\le i-1}B_\fr{\e}3(x_{i'k'})}f(u)+\fr{1}{2l+m}\int_{\O\sm\Cup_{i',k'}B_\fr{\e}3(x_{i'k'})}f(u)\q\q\q\q\q\q f(u)=\fr{e^u}{\int_\O e^u}.$$
To conclude the proof, we suffice to show that
$$\l|\int_\O(hf(u)-h\mrm d\mu(u))\r|\le\e\|h\|_{\mrm{Lip}\l(\ol\O\r)}.$$
We split the integral between the union of the balls of radius $\fr{\e}3$ and its complement: on the latter,
$$\l|\int_{\O\sm\Cup_{ik}B_\fr{\e}3(x_{ik})}(hf(u)-h\mrm d\mu(u))\r|=\l|\int_{\O\sm\Cup_{ik}B_\fr{\e}3(x_{ik})}hf(u)\r|\le\|h\|_{L^\infty(\O)}\int_{\O\sm\Cup_{ik}B_\fr{\e}3(x_{ik})}f(u)\le\fr{\e}3\|h\|_{\mrm{Lip}\l(\ol\O\r)}.$$
On the union of balls, we have:
\beqy
&&\l|\int_{\Cup_{ik}B_\fr{\e}3(x_{ik})}(hf(u)-h\mrm d\mu(u))\r|\\
&=&\l|\int_{\Cup_{ik}B_\fr{\e}3(x_{ik})}hf(u)-\sum_{ik}\l(\int_{B_\fr{\e}3(x_{ik})\sm\Cup_{k'\le k-1\tx{ or }i'\le i-1}B_\fr{\e}3(x_{i'k'})}f(u)+\fr{1}{2l+m}\int_{\O\sm\Cup_{i',k'}B_\fr{\e}3(x_{i'k'})}f(u)\r)h(x_{ik})\r|\\
&=&\l|\sum_{ik}\l(\int_{B_\fr{\e}3(x_{ik})\sm\Cup_{k'\le k-1\tx{ or }i'\le i-1}B_\fr{\e}3(x_{i'k'})}f(u)(h-h(x_{ik}))-\fr{h(x_{ik})}{2l+m}\int_{\O\sm\Cup_{i',k'}B_\fr{\e}3(x_{i'k'})}f(u)\r)\r|\\
&\le&\|\n h\|_{L^\infty(\O)}\sum_{ik}\int_{B_\fr{\e}3(x_{ik})\sm\Cup_{k'\le k-1\tx{ or }i'\le i-1}B_\fr{\e}3(x_{i'k'})}f(u)|\cd-x_{ik}|+\|h\|_{L^\infty(\O)}\int_{\O\sm\Cup_{i',k'}B_\fr{\e}3(x_{i'k'})}f(u)\\
&\le&\fr{\e}3\|\n h\|_{L^\infty(\O)}\int_{\Cup_{i',k'}B_\fr{\e}3(x_{i'k'})}f(u)+\fr{\e}3\|h\|_{L^\infty(\O)}\\
&\le&\fr{\e}3\|\n h\|_{L^\infty(\O)}+\fr{\e}3\|h\|_{L^\infty(\O)}\\
&\le&\fr{2}3\e\|h\|_{\mrm{Lip}'\l(\ol\O\r)}.
\eeqy
The proof is now complete.
\epf\

We are now in condition to prove the main result of this subsection.\\
We will construct $\Psi:\mcal J_{\b,\rho}^{-L}\to(\O_\pa)_{K,I}$ in the following way. The element in $(\O_\pa)_K$ will be given by the retraction $\wt\Psi$, while the element in $\S^{I-1}$ is just the normalization of $\Pi_Iu\in\R^I$. The choice of the third parameter in the join will be more delicate, especially in the homotopy, because we need to be sure that everything is well-defined outside the endpoints of the interval.

\bpf[Proof of Proposition \ref{psi}]${}$\\
Take $\e_0$ as in Lemma \ref{retr} and $L=L(\e_0)$ as in Lemma \ref{epsl}. We define the map $\Psi:\mcal J_{\b,\rho}^{-L}\to(\O_\pa)_{K,I}$ as
$$\Psi(u)=(\mu(u),\si(u),t(u)):=\l(\wt\Psi\l(\fr{e^u}{\int_\O e^u}\r),\fr{\Pi_Iu}{\|\Pi_Iu\|},\min\{1,\|\Pi_Iu\|\}\r).$$
We need to verify that it is well-posed, namely that $\mu(u)$ is well-defined if $t\ne1$ and $\si(u)$ is well-defined if $t\ne0$.\\
Assume $t\ne1$: this means $\|\Pi_Iu\|<1$ so, since $\mcal J_{\b,\rho}(u)\le-L$, Lemma \ref{deps} will give $d_{\mrm{Lip}'\l(\ol\O\r)}\l(\fr{e^u}{\int_\O e^u},(\O_\pa)_K\r)\le\e_0$; hence, Lemma \ref{retr} ensures that $\wt\Psi$ is well-defined, hence $\mu(u)$ is.\\
On the other hand, if $t\ne0$, then $\Pi_Iu\ne0$, hence one can define $\fr{\Pi_Iu}{\|\Pi_Iu\|}$.\\

As for second part of the lemma, consider the map $\Phi:=\Phi^{\La_0}$ as defined in Proposition \ref{test}, with $\La_0\gg1$ so large that $\Phi^{\La_0}((\O_\pa)_{K,I})\sub J_{\b,\rho}^{-L}$.\\
To get a homotopical equivalence, we let $\La$ go to $+\infty$. One immediately sees that $\fr{e^{\phi^{\La(1-t)}(\mu)}}{\int_\O e^{\phi^{\La(1-t)}(\mu)}}\us{\La\to+\infty}\wk\mu$ for any $\mu\in(\O_\pa)_K$ and $t\ne1$, hence being $e^{\psi^{\La t}}$ negligible with respect to $\int_\O e^{\phi^{\La(1-t)}}$ (see proof of Lemma \ref{exp}) one also has $\fr{e^{\Phi^\La(\z)}}{\int_\O e^{\Phi^\La(\z)}}\us{\La\to+\infty}\wk\mu$. Similarly, since $\phi^{\La(1-t)}-\fint_\O\phi^{\La(1-t)}$ is bounded in $L^2(\O)$, its projection will be negligible with respect to $\psi^{\La t}$, therefore $\Pi_I\Phi^\La\us{\La\to+\infty}\to\si$ as long as $t\ne0$.\\
The scalar parameter $t$ in the join will be more delicate to handle, because by the proof of Lemma \ref{sq} one gets $t\l(\Phi^\La(\z)\r)\sim\min\l\{1,\log^+(\La t)\r\}$; moreover, it is forced to be either $0$ or $1$ if either element in the join is not defined.\\
Therefore, before letting $\La$ go to $+\infty$, we need to properly rescale such a parameter, taking into account when it is allowed to be different from $0$ and/or $1$. To this purpose, we notice that, since $\fr{e^{\phi^{\La(1-t)}(\mu)}}{\int_\O e^{\phi^{\La(1-t)}(\mu)}}$ gets closer to $(\O_\pa)_{K,I}$ as $\La(1-t)$ is larger, we can assume that $\mu\l(\Phi^\La(\z)\r)$ is well-defined for $\La(1-t)\ge\fr{\La_0}3$ and similarly that $t\l(\Phi^\La(\z)\r)$ is well-defined for $\La t\ge\fr{\La_0}3$. Therefore, we will construct an intermediate parameter $t'$, which we set to be either $0$ or $1$ if $t$ is outside the previous range, which fills the whole $(0,1)$ as $\La$ goes to $+\infty$: as a first step, we fix $\La_0$ and interpolate linearly between $t\l(\Phi^\La(\z)\r)$ and $t'$, and everything is well defined because $\La=\La_0$ is fixed; then, we pass to the limit as $\La\to+\infty$ and, by the previous considerations, it is still well-posed and in the limit we recover the identity map.\\
Precisely, a continuous homotopical equivalence between $\Psi\c\Phi^{\La_0}$ and $\mrm{Id}_{(\O_\pa)_{K,I}}$ is given by
$$F(\z,s)=\l\{\bll\l(\mu\l(\Phi^{\La_0}(\z)\r),\si\l(\Phi^{\La_0}(\z)\r),(1-2s)t\l(\Phi^{\La_0}(\z)\r)+2st'(t,1)\r)&\tx{if }0\le s<\fr{1}2\\\l(\mu\l(\Phi^\fr{\La_0}{2-2s}(\z)\r),\si\l(\Phi^\fr{\La_0}{2-2s}(\z)\r),t'(t,2-2s)\r)&\tx{if }\fr{1}2\le s<1\\\z&\tx{if }s=1\earr\r.$$
$$\tx{with}\q\q\q t'(t,r):=\l\{\bll0&\tx{if }t<\fr{r}3\\\fr{3t-r}{3-2r}&\tx{if }\fr{r}3\le t\le1-\fr{r}3\\1&\tx{if }t>1-\fr{r}3\earr\r.$$
\epf\

\brem${}$\\
All the result shown in this section hold true, also when $K$ and/or $I$ equals zero. In each case, the space $(\O_\pa)_{K,I}$ is replaced by
$$X=\l\{\bll(\O_\pa)_K&\tx{if }I\ne K=0\\\S^{I-1}&\tx{if }K\ne I=0\\\es&\tx{if }I=K=0\earr\r..$$
One can easily see that all the proofs are still valid in all these cases. When $I=0$, in Proposition \ref{test} we just consider $\Phi^\La(\mu,-,0)$ and when $K=0$ we take $\Phi^\la(-,\si,1)$; in Proposition \ref{psi} we just set $F(\z,s)=\Phi\c\Phi^\fr{\La}{1-s}$.\\
If $I=K=0$, then Propositions \ref{test} and \ref{psi} make no sense, but Lemma \ref{mtimpr} applied with $l=0$, $m=1$, $\O_{21}=\O$ implies that $\mcal J_{\b,\rho}$ is coercive, namely $\mcal J_{\b,\rho}^{-L}=\es$ for $L$ large. For this reason, the proof of Theorem \ref{ex} can be adapted also to this case.
\erem\

\section{Proof of the main result}\

We need one last lemma concerning the Morse property of the functional $\mcal J_{\b,\rho}$.

\blem\label{morse}${}$\\
Assume \eqref{ks} has no non-trivial solutions and $\b-\fr{\rho}{|\O|}\ne-\la_j$ for any $j\in\N$.\\
Then $\mcal J_{\b,\rho}$ is a Morse functional and the Morse index $J$ of the trivial solution $u\eq0$ is such that $\la_{J+1}<\b-\fr{\rho}{|\O|}<\la_J$.
\elem\

\bpf${}$\\
One immediately sees that the second derivative of $\mcal J_{\b,\rho}$ is given by
$$\mcal J_{\b,\rho}(u)[v,w]=\int_\O(\n v\cd\n w)+\b\int_\O vw-\rho\frac{\int_\O vwe^u\int_\O e^u-\int_\O ve^u\int_\O we^u}{\l(\int_\O e^u\r)^2},$$
hence in $u\eq0$ its quadratic form is
$$\mcal J_{\b,\rho}(0)[v,v]=\int_\O|\n v|^2+\l(\b-\fr{\rho}{|\O|}\r)\int_\O v^2.$$
Assume the only solution to \eqref{ks} is the trivial one. Then, the $\mcal J_{\b,\rho}$ is a Morse functional if and only if the previous quadratic form is nondegenerate. One immediately sees that this depends on the relative position of $\fr{\rho}{|\O|}-\b$ and the $\la_j(\O)$'s as we required.
\epf\

Now we are finally in position to prove the main result of the paper.

\bpf[Proof of Theorem \ref{ex}]${}$\\
Assume, by contradiction, that $u\eq0$ is the only solution to \eqref{ks}.\\
By Lemma \ref{morse}, $\mcal J_{\b,\rho}$ is a Morse functional and the Morse index of the solution is $J$, therefore by Morse theory the relative homology of sublevels satisfies $H_q\l(\mcal J_{\b,\rho}^L,\mcal J_{\b,\rho}^{-L}\r)=\l\{\bll\Z&\tx{if }q=J\\0&\tx{if }q\ne J\earr\r.$, for any $L>0$.\\
Moreover, by Corollary \ref{subl}, $\mcal J_{\b,\rho}^L$ is contractible if $L$ is large enough; therefore, the exactness of the sequence (see \cite{hat}, Theorem $2.13$ and Proposition $2.22$)
$$\ds\to\wt H_q\l(\mcal J_{\b,\rho}^{-L}\r)\to\wt H_q\l(\mcal J_{\b,\rho}^L\r)\to H_q\l(\mcal J_{\b,\rho}^L,\mcal J_{\b,\rho}^{-L}\r)\to\wt H_{q-1}\l(\mcal J_{\b,\rho}^{-L}\r)\to\wt H_{q-1}\l(\mcal J_{\b,\rho}^L\r)\to\ds$$
yields
$$\wt H_q\l(\mcal J_{\b,\rho}^{-L}\r)=H_{q+1}\l(\mcal J_{\b,\rho}^L,\mcal J_{\b,\rho}^{-L}\r)=\l\{\bll\Z&\tx{if }q=J-1\\0&\tx{if }q\ne J-1\earr\r.,$$
which contradicts Proposition \ref{omo}. In fact, if $2K+I\ne J$, Corollary \ref{omosub} gives a non-trivial homology group for $q=2K+I-1\ne J-1$; moreover, if $\O$ is not simply connected and $K>0$, then even when $2K+I=J$ we get a bigger homology group: $\Z\sne\Z^{\bin{K+g}g}\inc H_{J-1}\l(\mcal J_{\b,\rho}^{-L}\r)$
\epf\

\section*{Acknowledgments}\

The author wishes to thank Professor Angela Pistoia for the discussions concerning the topics of the paper.

\bibliography{kellersegel}
\bibliographystyle{abbrv}

\end{document}